# An analytic bifurcation principle for Fredholm operators

Matthias Stiefenhofer*

ABSTRACT. Equations of the form $G[z] = 0$, $G: B \longrightarrow \bar{B}$ smooth, $B, \bar{B}$ real Banach spaces are investigated with the aim of continuing the basic solution $G[0] = 0$ with $G'[0]$ Fredholm operator to a solution curve $G[z(\varepsilon)] = 0$ with the implicit function theorem. If $G'[0]$ is surjective, then the transversality condition of the implicit function theorem can be satisfied in a straightforward way, yielding a regular solution curve, whereas otherwise the equation $G[z] = 0$ has to be extended appropriately for reaching a surjective linearization accessible to the implicit function theorem. This extension process, implying in the first step the standard bifurcation theorem of simple bifurcation points, is continued arbitrarily, yielding a sequence of bifurcation results presumably being applicable to bifurcation points with finite degeneracy.

**Keywords:** Implicit function theorem, Banach space, bifurcation theory, Fredholm operator

Contents

## *1. Introduction*

An equation of class $C^\infty$ of the form

$$G[z] = 0\,, \qquad G : B \longrightarrow \bar{B}\,, \qquad G[\bar{z}_0] = 0\,, \qquad G'[\bar{z}_0] \;\; \textit{Fredholm operator},$$

$B, \bar{B}$ real Banach spaces is considered, with the aim of finding solution curves $G[z(\varepsilon)] = 0, |\varepsilon| \ll 1$, through $\bar{z}_0$. Now, if there exists a solution curve $G[z(\varepsilon)] = 0$ with $z(0) = \bar{z}_0$, then the curve can be established locally out of the known solution $\bar{z}_0$ with the implicit function theorem under the assumption of

$$G'[\bar{z}_0] \in L[\,B, \bar{B}\,] \tag{1.1}$$

being surjective. Solution curves of this kind are typically called regular solution curves. On the other hand, if $G'[\bar{z}_0]$ is not surjective, then the existence of the curve can be shown provided that the extended linear sum operator

$$[\; G'[\bar{z}_0]\,,\; 2G''[\bar{z}_0] \cdot \bar{z}_1 \;] \;\in L[\; B \times N[\; G'[\bar{z}_0]\,]\,,\; \bar{B}\;]$$

*Seminar Bifurcation Theory, ETH Zürich, Valbella (1996)*

defined by

$$[\, G'[\bar{z}_0]\,,\;\; 2G''[\bar{z}_0]\cdot\bar{z}_1\,]\cdot\binom{b}{n} := G'[\bar{z}_0]\cdot b + 2G''[\bar{z}_0]\cdot\bar{z}_1 n \tag{1.2}$$

is surjective. Here, $\bar{z}_1$ denotes the tangential direction of $z(\varepsilon)$ in the central point $z(0) = \bar{z}_0$, whereas $N\big[G'[\bar{z}_0]\big]$ labels null space of $G'[\bar{z}_0]$. The regularity condition (1.2) occurs at simple bifurcation points.

The transition from condition (1.1) to condition (1.2) is accomplished by appending a second linear mapping, namely $2G''[\bar{z}_0]\cdot\bar{z}_1$, to the nonsurjective first mapping $G'[\bar{z}_0]$ which acts exactly on the subspace of $B$ not contributing to the range of $G'[\bar{z}_0]$ in $\bar{B}$, i.e. the second linear mapping acts on the kernel of the first linear mapping. In this paper, we show, how to continue this extension principle arbitrarily.

For this purpose and with the notation $\bar{z}_i := z^{(i)}(0), i \geq 0$, an iteratively defined chain of linear mappings

$$\begin{array}{llllllllll}
\tilde{S}_1(\bar{z}_0) & := & G'[\bar{z}_0] & \in & L[\,N_0,\bar{B}\,]\,, & N_0 & := & B & & \\
\tilde{S}_2(\bar{z}_1,\bar{z}_0) & := & 2G''[\bar{z}_0]\cdot\bar{z}_1 & \in & L[\,N_1,\bar{B}\,]\,, & N_1 & := & N[\,G'[\bar{z}_0]\,] & \subset & N_0 \\
\tilde{S}_3(\bar{z}_2,\bar{z}_1,\bar{z}_0) & := & \cdots & \in & L[\,N_2,\bar{B}\,]\,, & N_2 & := & \cdots & \subset & N_1 \;\subset\; N_0 \\
\vdots & & & & & \vdots & & & &
\end{array} \tag{1.3}$$

is constructed with associated sum operators satisfying the relation

$$R\big[\,\tilde{S}_1(\bar{z}_0)\,\big] \;\subset\; R\big[\;[\,\tilde{S}_1(\bar{z}_0)\,,\;\tilde{S}_2(\bar{z}_1,\bar{z}_0)\,]\;\big] \;\subset\; \cdots \;\subset\; \bar{B}\,.$$

Here $R[\cdot]$ denotes the range of the linear operator within the brackets. Now, it can be shown, if there exists a $k \geq 1$ with

$$R\big[\;[\,\tilde{S}_1(\bar{z}_0)\,,\;\cdots\,,\;\tilde{S}_{k+1}(\bar{z}_k,\ldots,\bar{z}_0)\,]\;\big] = \bar{B}\,, \tag{1.4}$$

then the curve $z(\varepsilon)$ can be established by implicit function theorem and the first $k$ derivatives $(\bar{z}_k,\ldots,\bar{z}_1)$ of $z(\varepsilon)$ in the known zero $\bar{z}_0$. In this sense, we say a solution curve $z(\varepsilon)$ is $k$-regular, if it satisfies the regularity condition (1.4).

Further, note that the derivatives $\bar{z}_1, \bar{z}_2, \ldots$ satisfy the following necessary conditions obtained by successively differentiating the implicit equation $G[z(\varepsilon)] = 0$

$$\begin{array}{llll}
T^0 & := & G[\bar{z}_0] & = \; 0 \\
T^1 & := & G'[\bar{z}_0]\cdot\bar{z}_1 & = \; 0 \\
T^2 & := & G''[\bar{z}_0]\cdot\bar{z}_1^2 \; + \; G'[\bar{z}_0]\cdot\bar{z}_2 & = \; 0 \\
\vdots & & &
\end{array}$$

whereas some sufficient conditions concerning the existence of solution curves of $G[z] = 0$ can be comprised in the following way:

> Each solution of $[T^0,\cdots,T^{2k}] = 0$, satisfying the regularity condition (1.4), ensures the existence of a nontrivial curve of solutions $z(\varepsilon)$ through $\bar{z}_0$. (1.5)

In case of $k = 1$, each solution of

$$[\,T^0,\;T^1,\;T^2\,] = [\,G[\bar{z}_0]\,,\;\;G'[\bar{z}_0]\cdot\bar{z}_1\,,\;\;G''[\bar{z}_0]\cdot\bar{z}_1^2 \;+\; G'[\bar{z}_0]\cdot\bar{z}_2\,] = 0\,,$$

satisfying (1.4), guarantees a nontrivial solution curve $z(\varepsilon)$. Thus, for $k = 1$, result (1.5) agrees with the classical bifurcation theorem of simple bifurcation points [2].

In case of $k = 2$, (1.5) corresponds to a bifurcation result in [3] dealing with bifurcation points of tangentially touching branches with different curvatures. Analogously, for arbitrary values of $k$, it is possible to prove the existence of curves agreeing in the first $k-1$ derivatives in the bifurcation point, but disagreeing in the $k$-th derivative. Moreover, the existence of solution curves with vanishing derivatives, e.g. cusp curves, can be shown.

In sections 2 and 3, the equations $T^i = 0, i \geq 1$, as well as the iteratively defined linear mappings $\tilde{S}_i(\bar{z}_{i-1}, \dots, \bar{z}_0)$ from (1.3) are constructed for finally proving (1.5) in section 4.

## *2. Results*

At first, the higher order chain rule [1] applied to $G[z(\varepsilon)] = 0$ reads for $k \geq 1$

$$T^k(z_k, \dots, z_0) = \sum_{\beta=1}^{k} G_0^\beta \sum_{\substack{n_1+\cdots+n_k=\beta \\ 1\cdot n_1+\cdots+k\cdot n_k=k}} \frac{k!}{n_1! \cdots n_k!} \prod_{\tau=1}^{k} \left(\frac{1}{\tau!} z_\tau\right)^{n_\tau} \in L[B, \bar{B}] \tag{2.1}$$

with $G_0^\beta$ denoting the $\beta$-th derivative of $G$ in $z_0$ and $T^k$ depending explicitly from $(z_k, \dots, z_0)$. However, in some situations $T^k$ will be also interpreted as a function depending on further dummy variables $z_{k+1}, z_{k+2}, \dots$ . Now, concerning the relevant system of equations from (1.5)

$$[\, T^0(z_0), \cdots, T^k(z_k, \dots, z_0),\ T^{k+1}(z_{k+1}, \dots, z_0), \cdots, T^{2k}(z_{2k}, \dots, z_0)\, ] = 0 \ \in\ \bar{B}^{2k+1}, \tag{2.2}$$

a direct inspection of (2.1) implies the following linear structure with respect to the last $k$ components $T^{k+1}(z_{k+1}, \dots, z_0), \dots, T^{2k}(z_{2k}, \dots, z_0)$

$$\begin{pmatrix} T^{2k}(z_{2k}, \dots, z_0) \\ \vdots \\ T^{k+1}(z_{k+1}, \dots, z_0) \end{pmatrix} = \Delta^k(z_k, \dots, z_0) \cdot \begin{pmatrix} z_{2k} \\ \vdots \\ z_{k+1} \end{pmatrix} + \mathrm{I}^k(z_k, \dots, z_0) \tag{2.3}$$

with an upper triangular matrix $\Delta^k(z_k, \dots, z_0) \in L[B^k, \bar{B}^k]$ and $\mathrm{I}^k(z_k, \dots, z_0) \in \bar{B}^k$ which are explicitly defined in the next section. Hence, there exists a solution $(\bar{z}_{2k}, \dots, \bar{z}_{k+1},\ \bar{z}_k, \dots, \bar{z}_0) \in B^{2k+1}$ of the original system (2.2) if and only if the conditions

$$[\, T^0(z_0), \dots, T^k(z_k, \dots, z_0)\, ] = 0 \ \in \bar{B}^{k+1} \quad \wedge \quad \mathrm{I}^k(z_k, \dots, z_0) \in R[\, \Delta^k(z_k, \dots, z_0)\, ] \tag{2.4}$$

can be satisfied for an element $(\bar{z}_k, \dots, \bar{z}_0) \in B^{k+1}$, suggesting a subordinate meaning of the remaining variables $(z_{2k}, \dots, z_{k+1}) \in B^k$ in (2.2), (2.3).

Now, if a solution $(\bar{z}_k, \dots, \bar{z}_0) \in B^{k+1}$ of (2.4) or equivalently a solution $(\bar{z}_{2k}, \dots, \bar{z}_{k+1},\ \bar{z}_k, \dots, \bar{z}_0) \in B^{2k+1}$ of (2.2) is known, then all solutions of equations (2.2) are obviously given by

$$[\, T^0, \dots, T^{2k}]\,[\, \begin{pmatrix} \bar{z}_{2k} \\ \vdots \\ \bar{z}_{k+1} \end{pmatrix} + \begin{pmatrix} b_{2k} \\ \vdots \\ b_{k+1} \end{pmatrix},\ \bar{z}_k, \dots, \bar{z}_0\ ] = 0\,, \quad \begin{pmatrix} b_{2k} \\ \vdots \\ b_{k+1} \end{pmatrix} \in N[\, \Delta^k(\bar{z}_k, \dots, \bar{z}_0)\, ]\,. \tag{2.5}$$

Using this affine subspace of solutions within $B^{2k+1}$ and the ansatz

$$G\left[\, z_0 + \varepsilon \cdot z_1 + \cdots + \tfrac{1}{(2k+1)!}\, \varepsilon^{2k+1} \cdot z_{2k+1} \right] = 0\,, \tag{2.6}$$

the existence of a solution curve $z(\varepsilon)$ with $z^{(i)}(0) = \bar{z}_i, i = 0, \dots, k$, is established in the following way. After a simple calculation, we see that the first $2k + 1$ coefficients of a Taylor expansion with respect to $\varepsilon$ of (2.6) agree with $T^0/0!, \dots, T^{2k}/(2k)!$, implying the remainder equation

$$T^{2k+1}[\, z_{2k+1}, \begin{pmatrix} \bar{z}_{2k} \\ \vdots \\ \bar{z}_{k+1} \end{pmatrix} + \begin{pmatrix} b_{2k} \\ \vdots \\ b_{k+1} \end{pmatrix}, \bar{z}_k, \dots, \bar{z}_0 \,] + \varepsilon \cdot r(\varepsilon, z_{2k+1}, b_{2k}, \dots, b_{k+1}) = 0 \qquad (2.7)$$

after inserting the solutions (2.5) into the Taylor expansion. Further, an inspection of the basic formula (2.1) shows the linearity property

$$T^{2k+1}(z_{2k+1}, \dots, z_0) = W^{2k+1}(z_k, \dots, z_0) \cdot \begin{pmatrix} z_{2k+1} \\ \vdots \\ z_{k+1} \end{pmatrix} + R^{2k+1}(z_k, \dots, z_0) \in \bar{B} \qquad (2.8)$$

concerning the leading term in (2.7) with $W^{2k+1}(z_k, \dots, z_0) \in L[B^{k+1}, \bar{B}]$ and $R^{2k+1}(z_k, \dots, z_0) \in \bar{B}$ which are again defined explicitly in the next section. Therefore, with $b_{2k+1} \coloneqq z_{2k+1} \in B$ and $(b_{2k}, \dots, b_{k+1}) \in N[\Delta^k(z_k, \dots, z_0)]$, the remainder equation (2.7) reads for $\varepsilon = 0$

$$T^{2k+1}[\, z_{2k+1}, \begin{pmatrix} \bar{z}_{2k} \\ \vdots \\ \bar{z}_{k+1} \end{pmatrix} + \begin{pmatrix} b_{2k} \\ \vdots \\ b_{k+1} \end{pmatrix}, \bar{z}_k, \dots, \bar{z}_0 \,]$$

$$= W^{2k+1}(\bar{z}_k, \dots, \bar{z}_0) \cdot [ \begin{pmatrix} 0 \\ \bar{z}_{2k} \\ \vdots \\ \bar{z}_{k+1} \end{pmatrix} + \begin{pmatrix} b_{2k+1} \\ b_{2k} \\ \vdots \\ b_{k+1} \end{pmatrix} ] + R^{2k+1}(\bar{z}_k, \dots, \bar{z}_0)$$

$$= W^{2k+1}(\bar{z}_k, \dots, \bar{z}_0) \cdot \begin{pmatrix} b_{2k+1} \\ \vdots \\ b_{k+1} \end{pmatrix} + T^{2k+1}(0, \bar{z}_{2k}, \dots, \bar{z}_0) = 0 \,.$$

Thus, if the linear mapping

$$W^{2k+1}(\bar{z}_k, \dots, \bar{z}_0) \in L[\, B \times N[\, \Delta^k(\bar{z}_k, \dots, \bar{z}_0) \,], \, \bar{B} \,]$$

is surjective and if additionally a closed subspace $N^c$ exists with

$$B \times N[\, \Delta^k(\bar{z}_k, \dots, \bar{z}_0) \,] = N^c \oplus N\left[ W^{2k+1}(\bar{z}_k, \dots, \bar{z}_0)_{|B \times N[\, \Delta^k(\bar{z}_k, \dots, \bar{z}_0) \,]} \right], \qquad (2.9)$$

then $W^{2k+1}(\bar{z}_k, \dots, \bar{z}_0) \in GL[N^c, \bar{B}]$, implying a unique solution of (2.7) in $N^c$ for $\varepsilon = 0$ that can uniquely be continued to $\varepsilon \neq 0$ by use of the implicit function theorem. In some more detail, there exists a locally defined function $(b_{2k+1}, \dots, b_{k+1})(\varepsilon) \in N^c, |\varepsilon| \ll 1$, of class $C^\infty$ with

$$G[\, \underbrace{\sum_{\mu=0}^{k} \tfrac{1}{\mu!} \varepsilon^\mu \cdot \bar{z}_\mu + \sum_{\mu=k+1}^{2k+1} \tfrac{1}{\mu!} \varepsilon^\mu \cdot \left( \bar{z}_\mu + b_\mu(\varepsilon) \right)}_{=: z(\varepsilon)} \,] = 0$$

and summarizing the following existence theorem is shown.

**Theorem 2.1 :** Every element $(\bar{z}_k, \dots, \bar{z}_0) \in B^{k+1}, k \geq 1$, satisfying

(i) $G'[\bar{z}_0]$ Fredholm operator

(ii) $[\, T^0(\bar{z}_0)\,, \cdots, T^k(\bar{z}_k, \dots, \bar{z}_0)\,] = 0 \in \bar{B}^{k+1} \quad \wedge \quad \mathrm{I}^k(\bar{z}_k, \dots, \bar{z}_0) \in R[\, \Delta^k(\bar{z}_k, \dots, \bar{z}_0)\,]$

(iii) $W^{2k+1}(\bar{z}_k, \dots, \bar{z}_0) \in L[\, B \times N[\, \Delta^k(\bar{z}_k, \dots, \bar{z}_0)\,],\ \bar{B}\,]\ surjective$ (2.10)

ensures a solution curve $z(\varepsilon)$ of class $C^\infty$ with

$$z^{(i)}(0) = \bar{z}_i\,,\ i = 0, \dots, k\,.$$

**Remark :** The existence of a closed subspace $N^c$ with (2.9) is a consequence of the Fredholm property of $G'[\bar{z}_0]$ as will be shown in section 4.

The above formulation of the sufficient conditions establishing a nontrivial solution curve $z(\varepsilon)$ does not show the iterative aspect of finding solution curves suggested in the introduction. Therefore, the following equivalent theorem is proved.

**Theorem 2.2 :** Every solution $(\bar{z}_{2k}, \dots, \bar{z}_0) \in B^{2k+1}, k \geq 1$, of

$$[\, T^0(\bar{z}_0)\,, \dots, T^{2k}(\bar{z}_{2k}, \dots, \bar{z}_0)\,] = 0 \in \bar{B}^{2k+1}$$

with $G'[\bar{z}_0]$ Fredholm operator gives rise to an iteratively defined sequence of $k+1$ linear mappings

$$\begin{array}{llll}
\tilde{S}_1(\bar{z}_0) & \in L[\, N_0, \bar{B}\,]\,, & N_0 := B & \\
\tilde{S}_2(\bar{z}_1, \bar{z}_0) & \in L[\, N_1, \bar{B}\,]\,, & N_1 := N[\, G'[\bar{z}_0]\,] & \subset N_0 \\
\vdots & \quad \vdots & \vdots & \quad \vdots \\
\tilde{S}_{k+1}(\bar{z}_k, \dots, \bar{z}_0) & \in L[\, N_k, \bar{B}\,]\,, & N_k := \quad \cdots & \subset N_{k-1} \subset \cdots \subset N_0
\end{array}$$

such that in case of

$$\left[\, \tilde{S}_1(\bar{z}_0)\,, \cdots, \tilde{S}_{k+1}(\bar{z}_k, \dots, \bar{z}_0)\,\right] \in L[\, N_0 \times \cdots \times N_k,\ \bar{B}\,]\ surjective, \qquad (2.11)$$

a solution curve $z(\varepsilon)$ of class $C^\infty$ exists with

$$z^{(i)}(0) = \bar{z}_i\,,\ i = 0, \dots, k\,.$$

The successive construction of the linear mappings is given in the next section. In [4], the structure of the proof of Theorem 2.2 is given and the iteration of section 3 is defined in a complete way.

## *3. Iteration*

In this section, the definitions of $\Delta^k(z_k, \dots, z_0) \in L[B^k, \bar{B}^k]$, $\mathrm{I}^k(z_k, \dots, z_0) \in \bar{B}^k$ from (2.3) and $W^{2k+1}(z_k, \dots, z_0) \in L[B^{k+1}, \bar{B}]$, $R^{2k+1}(z_k, \dots, z_0) \in \bar{B}$ from (2.8), as well as the iteratively defined linear mappings from Theorem 2.2 are stated. The basis for the definitions is given by (2.1), partly implying involved but constructive formulas.

First, the components of the upper triangular matrix $\Delta^k(z_k, \dots, z_0)$ and the associated inhomogeneity $\mathrm{I}^k(z_k, \dots, z_0)$ from (2.3) read

$$\Delta^k_{i,i}(z_k, \dots, z_0) := G^1_0 \ \in\ L[B, \bar{B}]\,,\ \ i = 1, \dots, k \qquad (3.1)$$

$$\Delta^k_{i,j}(z_k, \dots, z_0) := 0 \ \in\ L[B, \bar{B}]\,,\ \ i = 2, \dots, k, \quad j = 1, \dots, i-1$$

$$\Delta^k_{k-i+1,k-j+1}(z_k,\dots,z_0) := \frac{1}{(k+j)!}\sum_{\beta=1}^{k+i} G_0^\beta \sum_{\substack{n_1+\cdots+n_k+1=\beta \\ 1\cdot n_1+\cdots+k\cdot n_k+(k+j)\cdot 1=k+i}} \frac{k!}{n_1!\cdots n_k!}\prod_{\tau=1}^{k}\left(\frac{1}{\tau!}\,z_\tau\right)^{n_\tau}$$

$$\in\ L[B,\bar{B}]\,,\quad i=2,\dots,k,\quad j=1,\dots,i-1,$$

as well as

$$I^k_{k-i+1}(z_k,\dots,z_0) := \sum_{\beta=1}^{k+i} G_0^\beta \sum_{\substack{n_1+\cdots+n_k=\beta \\ 1\cdot n_1+\cdots+k\cdot n_k=k+i}} \frac{k!}{n_1!\cdots n_k!}\prod_{\tau=1}^{k}\left(\frac{1}{\tau!}\,z_\tau\right)^{n_\tau}\ \in\ \bar{B}\,,\ \ i=1,\dots,k\,.$$

Further, the $k+1$ components of the linear sum operator from (2.10)

$$W^{2k+1}(z_k,\dots,z_0) = \left[\,W^{2k+1}_{2k+1},\dots,W^{2k+1}_{k+1}\,\right](z_k,\dots,z_0)\ \in\ L[B^{k+1},\bar{B}] \tag{3.2}$$

assume the form

$$W^{2k+1}_\mu(z_k,\dots,z_0) := \frac{1}{\mu!}\sum_{\beta=1}^{2k+1} G_0^\beta \sum_{\substack{n_1+\cdots+n_k+1=\beta \\ 1\cdot n_1+\cdots+k\cdot n_k+\mu=2k+1}} \frac{(2k+1)!}{n_1!\cdots n_k!}\prod_{\tau=1}^{k}\left(\frac{1}{\tau!}\,z_\tau\right)^{n_\tau}\ \in\ L[B,\bar{B}] \tag{3.3}$$

for $\mu = 2k+1,\dots,k+1$, whereas the corresponding inhomogeneity from (2.8) reads

$$R^{2k+1}(z_k,\dots,z_0) \coloneqq \sum_{\beta=1}^{2k+1} G_0^\beta \sum_{\substack{n_1+\cdots+n_k=\beta \\ 1\cdot n_1+\cdots+k\cdot n_k=2k+1}} \frac{(2k+1)!}{n_1!\cdots n_k!}\prod_{\tau=1}^{k}\left(\frac{1}{\tau!}\,z_\tau\right)^{n_\tau}\ \in\ \bar{B}\,.$$

In addition, the even components of system (2.2) are given by

$$T^{2k}(z_k,\dots,z_0) = \left[\,W^{2k}_{2k},\dots,W^{2k}_{k}\,\right](z_{k-1},\dots,z_0)\cdot\begin{pmatrix} z_{2k} \\ \vdots \\ z_k\end{pmatrix} + R^{2k}(z_{k-1},\dots,z_0) + \frac{(2k)!}{2(k!)^2}\,G_0^2\cdot z_k^2\ \in\ \bar{B} \tag{3.4}$$

with

$$W^{2k}_\mu(z_{k-1},\dots,z_0) := \frac{1}{\mu!}\sum_{\beta=1}^{2k} G_0^\beta \sum_{\substack{n_1+\cdots+n_k+1=\beta \\ 1\cdot n_1+\cdots+(k-1)\cdot n_{k-1}+\mu=2k}} \frac{(2k)!}{n_1!\cdots n_{k-1}!}\prod_{\tau=1}^{k-1}\left(\frac{1}{\tau!}\,z_\tau\right)^{n_\tau}\ \in\ L[B,\bar{B}] \tag{3.5}$$

for $\mu = 2k,\dots,k$, as well as

$$R^{2k}(z_{k-1},\dots,z_0) \coloneqq \sum_{\beta=1}^{2k} G_0^\beta \sum_{\substack{n_1+\cdots+n_{k-1}=\beta \\ 1\cdot n_1+\cdots+(k-1)\cdot n_{k-1}=2k}} \frac{(2k)!}{n_1!\cdots n_{k-1}!}\prod_{\tau=1}^{k-1}\left(\frac{1}{\tau!}\,z_\tau\right)^{n_\tau}\ \in\ \bar{B}\,.$$

In particular, we obtain from (3.3) and (3.5) for $k \ge 1$

$$W^{2k+1}_{2k+1}(z_k,\dots,z_0) = W^{2k}_{2k}(z_{k-1},\dots,z_0) = G_0^1\ \in\ L[B,\bar{B}]\,. \tag{3.6}$$

Next, we start with the definition of the linear mappings from Theorem 2.2 according to

$$\bar{S}_1(\bar{z}_0) := \tilde{S}_1(\bar{z}_0) \coloneqq S_1(\bar{z}_0) := G_0^1 \in L[\,N_0,\ R_0^c\,]\quad with\quad N_0 \coloneqq B,\ \ R_0^c \coloneqq \bar{B}\,. \tag{3.7}$$

Due to the Fredholm property of $G_0^1$, two subspaces $N_1^c \subset N_0 = B$ and $R_1^c \subset R_0^c = \bar{B}$ can be chosen with direct sum decompositions

$$\begin{aligned} B = N_0 = N_1^c \quad &\oplus \; N[\, S_1(\bar{z}_0)\,] \; =: N_1^c \; \oplus \; N_1 \\ \bar{B} = R_0^c = R[\, S_1(\bar{z}_0)\,] \; \oplus \; R_1^c \quad & =: R_1 \; \oplus \; R_1^c \end{aligned} \tag{3.8}$$

and all subspaces closed with respect to the induced norms from $B$ and $\bar{B}$ with continuous projection operators

$$P_{R_1} \in L[\, \bar{B}, R_1\,] \qquad \textit{and} \qquad P_{R_1^c} = I_{\bar{B}} - P_{R_1} \in L[\, \bar{B}, R_1^c\,] \tag{3.9}$$

and $I_{\bar{B}}$ denoting identity in $\bar{B}$. Then, the operator $S_1(\bar{z}_0)$ satisfies by construction the regularity condition

$$S_1(\bar{z}_0) \in GL[\, N_1^c, R_1\,]\,. \tag{3.10}$$

In the next step, define the operators

$$\begin{aligned} \bar{S}_2(\bar{z}_1, \bar{z}_0) &\coloneqq 2G''[\bar{z}_0] \cdot \bar{z}_1 && \in L[\, B, \bar{B}\,] \\ \tilde{S}_2(\bar{z}_1, \bar{z}_0) &\coloneqq \bar{S}_2(\bar{z}_1, \bar{z}_0)_{|N_1} && \in L[\, N_1, \bar{B}\,] \\ S_2(\bar{z}_1, \bar{z}_0) &\coloneqq P_{R_1^c}\, \tilde{S}_2(\bar{z}_1, \bar{z}_0) && \in L[\, N_1, R_1^c\,] \end{aligned} \tag{3.11}$$

with corresponding decomposition derived from $S_2(\bar{z}_1, \bar{z}_0) \in L[\, N_1, R_1^c\,]$

$$\begin{aligned} N_1 = N_2^c \quad &\oplus \; N[\, S_2(\bar{z}_1, \bar{z}_0)\,] \; =: N_2^c \; \oplus \; N_2 \\ R_1^c = R[\, S_2(\bar{z}_1, \bar{z}_0)\,] \; \oplus \; R_2^c \quad & =: R_2 \; \oplus \; R_2^c \end{aligned} \tag{3.12}$$

satisfying the same properties as decomposition (3.8) with (3.9), (3.10) and index 1 replaced by index 2. Further, two triangular schemes of the form

$$\begin{array}{cccc} \boxed{d_{0,1}} & & & \\ d_{1,1} & & & \\ d_{2,1} & \boxed{d_{2,2}} & & \\ d_{3,1} & d_{3,2} & & \\ d_{4,1} & d_{4,2} & \boxed{d_{4,3}} & \\ \vdots & \vdots & \vdots & \\ l=1 & l=2 & l=3 & \cdots \end{array} \qquad\qquad \begin{array}{cccc} \boxed{c_{0,1}} & & & \\ c_{1,1} & & & \\ c_{2,1} & \boxed{c_{2,2}} & & \\ c_{3,1} & c_{3,2} & & \\ c_{4,1} & c_{4,2} & \boxed{c_{4,3}} & \\ \vdots & \vdots & \vdots & \\ l=1 & l=2 & l=3 & \cdots \end{array}$$

are defined column by column according to

$$\begin{aligned} l \geq 1 \quad &: \quad \boxed{d_{2l-2,l} \coloneqq \tfrac{2l-1}{l}} && \boxed{c_{2l-2,l} \coloneqq 1} \\ n \geq l \quad &: \quad d_{2n-1,l} \coloneqq \tfrac{2n}{2n+1-l} && c_{2n-1,l} \coloneqq c_{2n-2,l} \cdot d_{2n-2,l} \\ n \geq l \quad &: \quad d_{2n,l} \coloneqq \tfrac{2n+1}{2n+2-l} && c_{2n,l} \coloneqq c_{2n-1,l} \cdot d_{2n-1,l} \end{aligned} \tag{3.13}$$

with the rows of the schemes partly comprised to diagonal matrices

$$
\begin{aligned}
D^{2n-1} &:= Diag\left[\, d_{2n-1,1}\,,\ \cdots\ , d_{2n-1,n}\,\right] \ \in \mathbb{R}^{n,n} \\
D^{2n} &:= Diag\left[\, d_{2n,1}\,,\ \cdots\ , d_{2n,n}\,\right] \ \in \mathbb{R}^{n,n} \\
C^{2n-1} &:= Diag\left[\, c_{2n-1,1}\,,\ \cdots\ , c_{2n-1,n}\,\right] \ \in \mathbb{R}^{n,n} \\
C^{2n} &:= Diag\left[\, c_{2n,1}\,,\ \cdots\ , c_{2n,n}\,\right] \ \in \mathbb{R}^{n,n}\,.
\end{aligned} \tag{3.14}
$$

Then, using the abbreviations

$$M^3 := I_B \in L[\,B,B\,]$$

$$E^2 := \begin{pmatrix} E_{1,1} & E_{1,2} \\ E_{2,1} & E_{2,2} \end{pmatrix} := \begin{pmatrix} I_B & -S_1^{-1}(\bar{z}_0)\, P_{R_1}\, \bar{S}_2(\bar{z}_1,\bar{z}_0) \\ 0 & I_B \end{pmatrix} \in L[\,B^2,B^2\,] \tag{3.15}$$

$$\begin{pmatrix} a^3 & \bar{a}^3 \\ A^3 & \bar{A}^3 \end{pmatrix} := (C^3)^{-1} \cdot E^2 \cdot C^3 \ \in\ L[\,B^2,B^2\,]\,,$$

the iteration starting with $k = 2$ is well defined according to the following notation

$$
\begin{aligned}
\bar{S}_{k+1}(\bar{z}_k,\ldots,\bar{z}_0) &:= \left[\, W_{2k-1}^{2k},\ldots,W_k^{2k}\,\right](\bar{z}_{k-1},\ldots,\bar{z}_0) \cdot \begin{pmatrix} \bar{a}^{2k-1} \\ M^{2k-1} \cdot \bar{A}^{2k-1} \end{pmatrix} + \frac{(2k)!}{(k!)^2}\, G_0^2 \bar{z}_k \ \in\ L[\,B,\bar{B}\,] \\
\tilde{S}_{k+1}(\bar{z}_k,\ldots,\bar{z}_0) &:= \bar{S}_{k+1}(\bar{z}_k,\ldots,\bar{z}_0)_{|N_k} \ \in\ L[\,N_k,\bar{B}\,] \\
S_{k+1}(\bar{z}_k,\ldots,\bar{z}_0) &:= P_{R_k^c}\, \tilde{S}_{k+1}(\bar{z}_k,\ldots,\bar{z}_0) \ \in\ L[\,N_k,R_k^c\,]
\end{aligned} \tag{3.16}
$$

with corresponding decomposition

$$
\begin{aligned}
N_k &= N_{k+1}^c \oplus N[\,S_{k+1}(\bar{z}_k,\ldots,\bar{z}_0)\,] =: N_{k+1}^c \oplus N_{k+1} \\
R_k^c &= R[\,S_{k+1}(\bar{z}_k,\ldots,\bar{z}_0)\,] \oplus R_{k+1}^c =: R_{k+1} \oplus R_{k+1}^c
\end{aligned}
$$

satisfying again the same properties as decompositions (3.8), (3.12) with index 1,2 replaced by index $k+1$. Hence, setting $S_0 := 0$, the linear mappings $S_i, i = 1,\ldots,k+1$, are iteratively defined in such a way that the mapping $S_i$ acts between the kernel of the previous mapping $N[S_{i-1}] = N_{i-1}$ and a complement $R_{i-1}^c$ of the range $R[S_{i-1}] = R_{i-1}$ of this mapping.

Hence, by construction, we end up with the following decompositions of $B$ and $\bar{B}$ at iteration step $k \geq 1$.

$$
\begin{array}{ccccccccc}
B = & N_1^c & \oplus & N_2^c & \oplus \cdots \oplus & N_{k+1}^c & \oplus & N_{k+1} \\
 & \uparrow & & \uparrow & & \uparrow & & \\
 & \boxed{S_1} & & \boxed{S_2} & & \boxed{S_{k+1}} & & \\
 & \downarrow & & \downarrow & & \downarrow & & \\
\bar{B} = & R_1 & \oplus & R_2 & \oplus \cdots \oplus & R_{k+1} & \oplus & R_{k+1}^c
\end{array}
$$

$$
N_0 = N_1^c \oplus N_1, \quad N_1 = N_2^c \oplus \cdots \oplus N_{k+1}, \quad \ldots, \quad N_k = N_{k+1}^c \oplus N_{k+1}
$$

$$
R_0^c = R_1 \oplus \cdots \oplus R_{k+1}^c, \quad R_1^c = R_2 \oplus \cdots \oplus R_{k+1}^c, \quad \ldots, \quad R_k^c = R_{k+1} \oplus R_{k+1}^c
$$

Figure 1 : Direct sum decompositions of $B$ and $\bar{B}$ at iteration step $k \geq 1$.

Note also that the relations

$$
S_i := P_{R_{i-1}^c}\, \bar{S}_{i|N_{i-1}} \in L[\,N_{i-1}, R_{i-1}^c\,] \qquad and \qquad S_i \in GL[\,N_i^c, R_i\,]\,, \;\; i = 1, \ldots, k+1
$$

are valid, where the bijectivity of $S_i, i = 1, \ldots, k+1$, between the subspaces $N_i^c$ and $R_i$ is depicted by arrows in figure 1.

Now, for closing iteration (3.16), the formulas concerning $\bar{a}^{2k+1}$, $\bar{A}^{2k+1}$ and $M^{2k+1}$ have to be established according to

$$
\begin{aligned}
& E_{k+1,k+1} := I_B \in L[\,B, B\,] \\
& E_{i,k+1} := -S_i^{-1}\, P_{R_i} \sum_{\nu=i+1}^{k+1} \bar{S}_\nu\, E_{\nu,k+1} \in L[\,B, B\,]\,, \;\; i = k, \ldots, 1 \\
& E^{k+1} := \begin{pmatrix} E_{1,1} & \cdots & E_{1,k+1} \\ & \ddots & \vdots \\ & & E_{k+1,k+1} \end{pmatrix} \in L[\,B^{k+1}, B^{k+1}\,] \\
& \begin{pmatrix} a^{2k} & \bar{a}^{2k} \\ A^{2k} & \bar{A}^{2k} \end{pmatrix} := \left(C^{2k}\right)^{-1} \cdot E^k \cdot C^{2k} \;\in\; L[\,B^k, B^k\,] \\
& \begin{pmatrix} a^{2k+1} & \bar{a}^{2k+1} \\ A^{2k+1} & \bar{A}^{2k+1} \end{pmatrix} := \left(C^{2k+1}\right)^{-1} \cdot E^{k+1} \cdot C^{2k+1} \;\in\; L[\,B^{k+1}, B^{k+1}\,]
\end{aligned}
\tag{3.17}
$$

with

$$
a^{2k} \in L[\,B^{k-1}, B\,]\,, \quad \bar{a}^{2k} \in L[\,B, B\,]\,, \quad A^{2k} \in L[\,B^{k-1}, B^{k-1}\,]\,, \quad \bar{A}^{2k} \in L[\,B, B^{k-1}\,]
$$

$$
a^{2k+1} \in L[\,B^k, B\,]\,, \quad \bar{a}^{2k+1} \in L[\,B, B\,]\,, \quad A^{2k+1} \in L[\,B^k, B^k\,]\,, \quad \bar{A}^{2k+1} \in L[\,B, B^k\,]
$$

and

$$M^{2k+1} := \begin{pmatrix} & & (a^{2k} \quad \bar{a}^{2k}) \\ & a^{2k-1} \cdot & (A^{2k} \quad \bar{A}^{2k}) \\ M^{2k-1} \cdot & A^{2k-1} \cdot & (A^{2k} \quad \bar{A}^{2k}) \end{pmatrix}_{\overline{k+1}} \in L[\,B^k, B^k\,]\,, \tag{3.18}$$

where the index $\overline{k+1}$ denotes cancelling the last row of the $(k+1) \times k$ matrix in brackets. Then, the iteration (3.16)-(3.18) can be continued with (3.16) and $k$ replaced by $k+1$.

Finally, a simple calculation shows the upper triangular property

$$M^{2k+1}_{i,j} = 0\,, \quad i = 2, \dots, k\,, \quad j = 1, \dots, i-1 \qquad and \qquad M^{2k+1}_{i,i} = I_B\,, \quad i = 1, \dots, k \tag{3.19}$$

implying

$$M^{2k+1} \in GL[\,B^k, B^k\,]\,. \tag{3.20}$$

**Remark :** The iteratively defined mappings $E_{i,k+1}$ in (3.17) can also be written in an explicit way according to

$$E_{k+1,k+1} = I_B\,, \qquad E_{k,k+1} = -S_k^{-1}\, P_{R_k}\, \bar{S}_{k+1}$$

$$E_{i,k+1} := -S_i^{-1}\, P_{R_i} \left( I_B + \sum_{\nu=1}^{k-i} (-1)^\nu \sum_{i+1 \le n_1 < \cdots < n_\nu \le k} \; \prod_{\tau=1}^{\nu} \bar{S}_{n_\tau}\, S_{n_\tau}^{-1}\, P_{R_{n_\tau}} \right) \bar{S}_{k+1}\,, \quad i = k-1, \dots, 1.$$

## *4. Proofs*

In the following Lemma, the main ingredients for proving Theorem 2.1 and Theorem 2.2 from section 2 are comprised.

**Lemma 4.1 :** Assume an element $(\bar{z}_{2k}, \dots, \bar{z}_0) \in B^{2k+1}$ with $[\,T^{2k}, \dots, T^0\,](\bar{z}_{2k}, \dots, \bar{z}_0) = 0$ and $k \ge 1$. Then, we obtain

(i) $$[\,T^{2k}, \dots, T^{k+1}\,](z_{2k}, \dots, z_{k+1}, \bar{z}_k, \dots, \bar{z}_0) = 0$$

$$\Leftrightarrow \begin{pmatrix} z_{2k} \\ \vdots \\ z_{k+1} \end{pmatrix} = \begin{pmatrix} \bar{z}_{2k} \\ \vdots \\ \bar{z}_{k+1} \end{pmatrix} + M^{2k+1} \cdot \begin{pmatrix} n_1 \\ \vdots \\ n_k \end{pmatrix}, \quad \begin{pmatrix} n_1 \\ \vdots \\ n_k \end{pmatrix} \in N_1 \times \cdots \times N_k$$

(ii) $$W^{2k+1}(\bar{z}_k, \dots, \bar{z}_0) \cdot \begin{pmatrix} I_B & \\ & M^{2k+1} \end{pmatrix} = [\,\bar{S}_1(\bar{z}_0)\,, \; \cdots \;, \; \bar{S}_{k+1}(\bar{z}_k, \dots, \bar{z}_0)\,] \cdot C^{2k+1}$$

(iii) $$N[\,[\,\bar{S}_1(\bar{z}_0)\,, \; \cdots \;, \; \bar{S}_{k+1}(\bar{z}_k, \dots, \bar{z}_0)\,]_{|N_0 \times \cdots \times N_k}\,] = R[\,E^{k+1}{}_{|N_1 \times \cdots \times N_{k+1}}\,]$$

(iv) $$M^{2k+1} \cdot \begin{pmatrix} d_{2k+1,2} & & \\ & \ddots & \\ & & d_{2k+1,k+1} \end{pmatrix} = \begin{pmatrix} d_{2k+1,2} & & \\ & \ddots & \\ & & d_{2k+1,k+1} \end{pmatrix} \cdot \begin{pmatrix} & a^{2k+1} \\ M^{2k+1} \cdot & A^{2k+1} \end{pmatrix}_{\overline{k+1}}$$

Note that the assumption $[\,T^{2k}, \dots, T^0\,](\bar{z}_{2k}, \dots, \bar{z}_0) = 0$ implies $[\,T^{2m}, \dots, T^0\,](\bar{z}_{2m}, \dots, \bar{z}_0) = 0$ for $m = 1, \dots, k$. Hence, if one of the assertions (i)-(iv) is true with $k$, then it is also true with $k$ replaced by $m = 1, \dots, k$.

Equivalence (i) states that if a solution $(z_{2k},\dots,z_{k+1},\ \bar{z}_k,\dots,\bar{z}_0)$ of $[T^{2k},\dots,T^0]=0$ is known, then all solutions of $[T^{2k},\dots,T^0]=0$ can be calculated, provided the leading coefficients $(\bar{z}_k,\dots,\bar{z}_0)$ remain fixed. These solutions are given by an affine subspace within $B^k$, defined by higher coefficients $(\bar{z}_{2k},\dots,\bar{z}_{k+1})$ of the assumed solution and parametrized by kernels $N_1,\dots,N_k$ of the linear operators $S_1(\bar{z}_0),\dots,S_k(\bar{z}_{k-1},\dots,\bar{z}_0)$.

Secondly, (ii) ascertains the link between the global approach of Theorem 2.1 and the viewpoint taken in Theorem 2.2 of iteratively constructing linear $S$ mappings until surjectivity of the sum operator is reached.

The statements (iii) and (iv) are of more technical nature, intensively employing the definitions from section 3.

**Proof of Theorems 2.1 and 2.2 :** First, we finish the theorems in section 2 by a few simple calculations using Lemma 4.1. First, (i) and (2.3) imply

$$N[\,\Delta^k(\bar{z}_k,\dots,\bar{z}_0)\,]=R[\,M^{2k+1}{}_{|N_1\times\cdots\times N_k}\,], \tag{4.1}$$

proving by (ii) and (3.16) the equivalence of the surjectivity conditions (2.10) and (2.11) of Theorem 2.1 and Theorem 2.2 respectively. Hence, the theorems are shown, if a closed subspace $N^c$ with (2.9) exists, where a sufficient condition reads

$$dim\left\{\ N\left[\,W^{2k+1}(\bar{z}_k,\dots,\bar{z}_0)_{|B\times N[\,\Delta^k(\bar{z}_k,\dots,\bar{z}_0)\,]}\right]\right\}<\infty\,. \tag{4.2}$$

Now, from $G_0^1$ Fredholm operator and $N_1=N[G_0^1]\supset N_2\supset\cdots\supset N_k$, we conclude $dim\{N_1\times\cdots\times N_k\}<\infty$ as well as

$$dim\{\,N[\,\Delta^k(\bar{z}_k,\dots,\bar{z}_0)\,]\,\}=dim\{\,R[\,M^{2k+1}{}_{|N_1\times\cdots\times N_k}\,]\,\}<\infty$$

by (4.1), yielding with respect to the null space of $W^{2k+1}$ in (4.2) under consideration of (3.2)

$$\underbrace{W^{2k+1}_{2k+1}}_{\substack{=G_0^1\\(3.6)}}\cdot b\;+\;\underbrace{[\,W^{2k+1}_{2k},\dots,W^{2k+1}_{k+1}\,]\cdot M^{2k+1}}_{=:\,A}\cdot\underbrace{\begin{pmatrix}n_1\\ \vdots\\ n_k\end{pmatrix}}_{=:\,n}$$

$$=\;G_0^1\cdot b\;+\;A\cdot n=[\,G_0^1\,,\;\;A\,]\cdot\begin{pmatrix}b\\ n\end{pmatrix}=0$$

and $b\in B,\ n\in N:=N_1\times\cdots\times N_k\,,\ A\in L[N,\bar{B}]$. Next, choose decompositions according to

$$B=N_1^c\quad\oplus\overbrace{N[\,G_0^1\,]}^{dim<\infty}\qquad\qquad N=\overbrace{N_A^c}^{dim<\infty}\ \oplus\ N[\,A\,]$$

$$\bar{B}=R[\,G_0^1\,]\ \oplus\underbrace{R_1^c}_{dim<\infty}\qquad\qquad \bar{B}=\underbrace{R[\,A\,]}_{dim<\infty}\ \oplus\ R_A^c$$

and note $G_0^1\in GL[N_1^c,R[G_0^1]]$ and $A\in GL\big[N_A^c,R[A]\big]$. Then, with subspace $U:=R[G_0^1]\cap R[A]\subset\bar{B}$, we obtain

$$[\,G_0^1\,,\;\;A\,]\cdot\begin{pmatrix}b\\ n\end{pmatrix}=0\quad\Leftrightarrow\quad\begin{pmatrix}b\\ n\end{pmatrix}=\begin{pmatrix}(G_0^1)^{-1}\\ -A^{-1}\end{pmatrix}\cdot u+\begin{pmatrix}n_1\\ n_A\end{pmatrix}$$

with $u\in U,\ n_1\in N[G_0^1],\ n_A\in N[A]$ and all subspaces of finite dimension, implying (4.2) as desired.

The Fredholm property of $G_0^1$, as well as the finite dimensionality of the subspaces is not necessary for the procedure to work. It is enough to require the closure of the subspaces within the decompositions. This aspect is treated in some more detail in [4].

In figure 2 an overview concerning the proof of Lemma 4.1 is given.

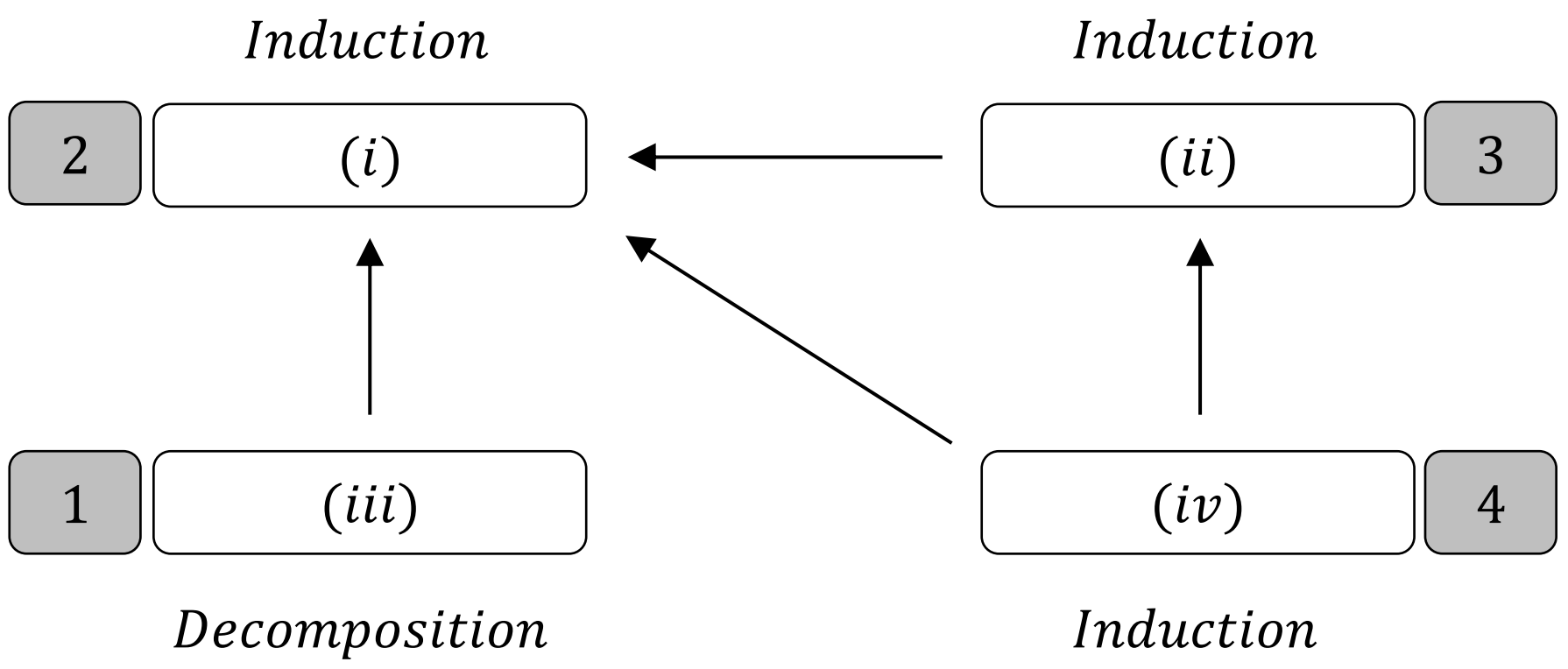

Figure 2 : The strategy for proving Lemma 4.1

First, we obtain (iii) by direct calculation using the decomposition shown in figure 1 and some essential definitions of the iteration process within section 3. Next, (i) is obtained by an induction argument relying on the validity of (ii), (iii) and (iv). Thirdly, (ii) is proved, again by induction as well as using (iv) as principal argument. Finally, the most technical part (iv) is shown by induction.

**Proof of Lemma 4.1 (iii) :** Dropping arguments, we obtain by the decomposition of figure 1

$$[\,\bar{S}_1, \cdots, \bar{S}_{k+1}\,] \cdot \begin{pmatrix} n_0 \\ \vdots \\ n_k \end{pmatrix} = 0$$

$$\Leftrightarrow \quad P_{R_i}[\,\bar{S}_1, \cdots, \bar{S}_{k+1}\,] \cdot \begin{pmatrix} n_0 \\ \vdots \\ n_k \end{pmatrix} = 0\,,\, i = 1, \dots, k+1 \quad \wedge \quad P_{R_{k+1}^c}[\,\bar{S}_1, \cdots, \bar{S}_{k+1}\,] \cdot \begin{pmatrix} n_0 \\ \vdots \\ n_k \end{pmatrix} = 0$$

$$\Leftrightarrow \quad \begin{pmatrix} \boxed{S_1} & P_{R_1}\bar{S}_2 & \cdots & P_{R_1}\bar{S}_{k+1} \\ & \ddots & \ddots & \vdots \\ & & \boxed{S_k} & P_{R_k}\bar{S}_{k+1} \\ & & & \boxed{S_{k+1}} \end{pmatrix} \cdot \begin{pmatrix} n_0 \\ \vdots \\ n_{k-1} \\ n_k \end{pmatrix} = 0$$

$$\wedge \quad \bar{S}_1 n_0 + \cdots + \bar{S}_{k+1} n_k \in R_1 \oplus \cdots \oplus R_{k+1}\,,$$

where the definitions of linear $S$ mappings in (3.16) imply $\bar{S}_1 n_0 \in R_1, \dots, \bar{S}_{k+1} n_k \in R_1 \oplus \cdots \oplus R_{k+1}$ and we can restrict to the first condition. Then, by the definitions of the operators $E_{i,j}$ in (3.17), the following equivalences result from bottom up solution of the triangular system

$$\Leftrightarrow \quad \begin{cases} n_k = \bar{n}_{k+1} = E_{k+1,k+1} \cdot \bar{n}_{k+1}\,, & \bar{n}_{k+1} \in N_{k+1} \\ n_{k-1} = \bar{n}_k - S_k^{-1} P_{R_k} \bar{S}_{k+1} \cdot \bar{n}_{k+1} = \left[\, E_{k,k} \quad E_{k,k+1} \,\right] \cdot \begin{pmatrix} \bar{n}_k \\ \bar{n}_{k+1} \end{pmatrix}, & \bar{n}_k \in N_k \\ \vdots \qquad \vdots \qquad \vdots & \vdots \end{cases}$$

$$\Leftrightarrow \quad \begin{pmatrix} n_0 \\ \vdots \\ n_k \end{pmatrix} = \begin{pmatrix} E_{1,1} & \cdots & E_{1,k+1} \\ & \ddots & \vdots \\ & & E_{k+1,k+1} \end{pmatrix} \cdot \begin{pmatrix} \bar{n}_1 \\ \vdots \\ \bar{n}_{k+1} \end{pmatrix}, \quad \begin{pmatrix} \bar{n}_1 \\ \vdots \\ \bar{n}_{k+1} \end{pmatrix} \in N_1 \times \cdots \times N_{k+1}$$

with the last equivalence inferred from a simple induction argument, thus finishing the proof of Lemma 4.1 (iii).

For later use, some easy to show preliminaries are added, based on (3.3), (3.5) and (3.13)

$$d_{m,1} = c_{m,1} = 1\,, \quad m \geq 0 \tag{4.3}$$

$$\frac{d_{2m,2+j}}{d_{2m-1,1+j}} = d_{2m,2}\,, \quad \frac{d_{2m+1,2+j}}{d_{2m,1+j}} = d_{2m+1,2}\,, \quad m \geq 1\,, \quad j = 0, \dots, m-1 \tag{4.4}$$

$$\frac{d_{2m+1,3+j}}{d_{2m-1,1+j}} = d_{2m,2} \cdot d_{2m+1,2}\,, \quad m \geq 2\,, \quad j = 0, \dots, m-2 \tag{4.5}$$

$$[\, W^{2m+1}_{2m+1}, \dots, W^{2m+1}_{m+2}\,](z_m, \dots, z_0) = [\, W^{2m}_{2m}, \dots, W^{2m}_{m+1}\,](z_{m-1}, \dots, z_0) \cdot D^{2m}\,, \quad m \geq 2 \tag{4.6}$$

$$W^{2m+1}_{m+1}(z_m, \dots, z_0) = [\, W^{2m}_{m}(z_{m-1}, \dots, z_0) + \frac{(2m)!}{(m!)^2}\, G^2_0 \cdot z_m\,] \cdot d_{2m,m+1}\,, \quad m \geq 2 \tag{4.7}$$

$$[\, W^{2m}_{2m}, \dots, W^{2m}_{m+1}\,](z_{m-1}, \dots, z_0) = [\, W^{2m-1}_{2m-1}, \dots, W^{2m-1}_{m}\,](z_{m-1}, \dots, z_0) \cdot D^{2m-1}, \quad m \geq 2 \tag{4.8}$$

**Proof of Lemma 4.1 (i) :** Next, according to figure 2, (i) is derived from (ii), (iii) and (iv) by induction from $m = 1$ to $m = k$. For $m = 1$, the equivalence from (i) is obviously true by (3.15) and

$$T^2(z_2, \bar{z}_1, \bar{z}_0) = G^2_0 \cdot \bar{z}^2_1 + G^1_0 \cdot z_2 = 0 \quad \Leftrightarrow \quad z_2 = \bar{z}_2 + I_B \cdot n_1 = \bar{z}_2 + M^3 \cdot n_1\,.$$

Now, suppose the equivalence (i) with $k$ replaced by $m - 1 \in \{1, \dots, k-2\}$

$$[\, T^{2m-2}, \dots, T^m\,](z_{2m-2}, \dots, z_m,\ \bar{z}_{m-1}, \dots, \bar{z}_0) = 0 \tag{4.9}$$

$$\Leftrightarrow \quad \begin{pmatrix} z_{2m-2} \\ \vdots \\ z_m \end{pmatrix} = \begin{pmatrix} \bar{z}_{2m-2} \\ \vdots \\ \bar{z}_m \end{pmatrix} + M^{2m-1} \cdot \begin{pmatrix} n_1 \\ \vdots \\ n_{m-1} \end{pmatrix}, \quad \begin{pmatrix} n_1 \\ \vdots \\ n_{m-1} \end{pmatrix} \in N_1 \times \cdots \times N_{m-1}\,.$$

Then, this affine subspace of solutions is successively plugged into the equations $T^{2m-1} = 0$ and $T^{2m} = 0$ for showing (i) with $k$ replaced by $m$. First, denoting $z_{2m-1} = n_0 \in N_0 = B$, the next equation reads

$$T^{2m-1}\,[\, n_0, \begin{pmatrix} \bar{z}_{2m-2} \\ \vdots \\ \bar{z}_m \end{pmatrix} + M^{2m-1} \cdot \begin{pmatrix} n_1 \\ \vdots \\ n_{m-1} \end{pmatrix}, \bar{z}_{m-1}, \dots, \bar{z}_0\,]$$

$$\overset{(2.8)}{=} W^{2m-1}(\bar{z}_{m-1}, \dots, \bar{z}_0) \cdot [\, \begin{pmatrix} 0 \\ \bar{z}_{2m-2} \\ \vdots \\ \bar{z}_m \end{pmatrix} + \begin{pmatrix} I_B & \\ & M^{2m-1} \end{pmatrix} \cdot \begin{pmatrix} n_0 \\ n_1 \\ \vdots \\ n_{m-1} \end{pmatrix}\,] + R^{2m-1}(\bar{z}_{m-1}, \dots, \bar{z}_0)$$

$$\overset{(2.8)}{=} W^{2m-1}(\bar{z}_{m-1}, \dots, \bar{z}_0) \cdot \begin{pmatrix} I_B & \\ & M^{2m-1} \end{pmatrix} \cdot \begin{pmatrix} n_0 \\ \vdots \\ n_{m-1} \end{pmatrix} + T^{2m-1}[\, 0,\, \bar{z}_{2m-2}, \dots, \bar{z}_0 \,]$$

$$\overset{(ii)}{=} [\, \bar{S}_1(\bar{z}_0), \cdots, \bar{S}_m(\bar{z}_{m-1}, \dots, \bar{z}_0) \,] \cdot C^{2m-1} \cdot \begin{pmatrix} n_0 \\ \vdots \\ n_{m-1} \end{pmatrix} - \bar{S}_1(\bar{z}_0) \cdot \bar{z}_{2m-1} = 0 \qquad (4.10)$$

under consideration of

$$W^{2m-1}_{2m-1}(z_{m-1}, \dots, z_0) = G^1_0 = \bar{S}_1(\bar{z}_0)$$

according to (3.6), (3.7) and

$$0 = T^{2m-1}[\, \bar{z}_{2m-1}, \bar{z}_{2m-2}, \dots, \bar{z}_0 \,] = \bar{S}_1(\bar{z}_0) \cdot \bar{z}_{2m-1} + T^{2m-1}[\, 0,\, \bar{z}_{2m-2}, \dots, \bar{z}_0 \,]\,.$$

Thus, by (iii) and the definition of $E^m$ in (3.15), (3.17), the affine subspace of solutions of (4.10) is uniquely given by

$$\begin{pmatrix} n_0 \\ n_1 \\ \vdots \\ n_{m-1} \end{pmatrix} = \begin{pmatrix} \bar{z}_{2m-1} \\ 0 \\ \vdots \\ 0 \end{pmatrix} + \begin{pmatrix} a^{2m-1} & \bar{a}^{2m-1} \\ A^{2m-1} & \bar{A}^{2m-1} \end{pmatrix} \cdot \begin{pmatrix} \bar{n}_1 \\ \bar{n}_2 \\ \vdots \\ \bar{n}_m \end{pmatrix} \qquad (4.11)$$

with $(\bar{n}_1, \dots, \bar{n}_m) \in N_1 \times \cdots \times N_m$ and $(\bar{z}_{2m-1}, 0, \dots, 0)^T$ representing a basic solution of (4.10).

Summing up (4.9) and (4.11), the solutions of the extended system

$$[\, T^{2m-1}, T^{2m-2}, \dots, T^{m+1}, T^m \,](z_{2m-1}, z_{2m-2}, \dots, z_m, \bar{z}_{m-1}, \dots, \bar{z}_0) = 0 \qquad (4.12)$$

are exactly given by

$$\begin{pmatrix} z_{2m-1} \\ z_{2m-2} \\ \vdots \\ z_m \end{pmatrix} = \begin{pmatrix} \bar{z}_{2m-1} + (a^{2m-1} \quad \bar{a}^{2m-1}) \cdot \begin{pmatrix} \bar{n}_1 \\ \vdots \\ \bar{n}_m \end{pmatrix} \\ \begin{pmatrix} \bar{z}_{2m-2} \\ \vdots \\ \bar{z}_m \end{pmatrix} + M^{2m-1} \cdot (A^{2m-1} \quad \bar{A}^{2m-1}) \cdot \begin{pmatrix} \bar{n}_1 \\ \vdots \\ \bar{n}_m \end{pmatrix} \end{pmatrix}$$

$$= \begin{pmatrix} \bar{z}_{2m-1} \\ \vdots \\ \bar{z}_m \end{pmatrix} + \begin{pmatrix} (a^{2m-1} \quad \bar{a}^{2m-1}) \\ M^{2m-1} \cdot (A^{2m-1} \quad \bar{A}^{2m-1}) \end{pmatrix} \cdot \begin{pmatrix} \bar{n}_1 \\ \vdots \\ \bar{n}_m \end{pmatrix}$$

with $(\bar{n}_1, \dots, \bar{n}_m) \in N_1 \times \cdots \times N_m$. Further, by (3.17) and (3.19), the last component reads $z_m = \bar{z}_m + \bar{n}_m$, $\bar{n}_m \in N_m$ and choosing $\bar{n}_m = 0$ as well as dropping the equation $T^m = 0$ in (4.12), the zeros of

$$[\, T^{2m-1}, T^{2m-2}, \dots, T^{m+1} \,](z_{2m-1}, z_{2m-2}, \dots, z_{m+1}, \bar{z}_m, \bar{z}_{m-1}, \dots, \bar{z}_0) = 0 \qquad (4.13)$$

are precisely given by

$$\begin{pmatrix} z_{2m-1} \\ \vdots \\ z_{m+1} \end{pmatrix} = \begin{pmatrix} \bar{z}_{2m-1} \\ \vdots \\ \bar{z}_{m+1} \end{pmatrix} + \begin{pmatrix} a^{2m-1} \\ M^{2m-1} \cdot A^{2m-1} \end{pmatrix}_{\bar{m}} \cdot \begin{pmatrix} \bar{n}_1 \\ \vdots \\ \bar{n}_{m-1} \end{pmatrix} \qquad (4.14)$$

with $(\bar{n}_1, \dots, \bar{n}_{m-1}) \in N_1 \times \cdots \times N_{m-1}$ and index $\overline{m}$ cancelling last row as usual.

In the next step, the solutions from (4.14) are plugged into the equation $T^{2m} = 0$, implying

$$0 = T^{2m}\Big[\underbrace{z_{2m}}_{=n_0 \in B}, \begin{pmatrix} \bar{z}_{2m-1} \\ \vdots \\ \bar{z}_{m+1} \end{pmatrix} + \begin{pmatrix} & a^{2m-1} \\ M^{2m-1} \cdot A^{2m-1} & \end{pmatrix}_{\overline{m}} \cdot \begin{pmatrix} \bar{n}_1 \\ \vdots \\ \bar{n}_{m-1} \end{pmatrix}, \bar{z}_m, \dots, \bar{z}_0\Big] \tag{4.15}$$

$$\overset{(3.4)}{\underset{(3.6)}{=}} \bar{S}_1(\bar{z}_0) \cdot n_0 + [W^{2m}_{2m-1}, \dots, W^{2m}_{m+1}](\bar{z}_{m-1}, \dots, \bar{z}_0) \cdot \begin{pmatrix} & a^{2m-1} \\ M^{2m-1} \cdot A^{2m-1} & \end{pmatrix}_{\overline{m}} \cdot \begin{pmatrix} \bar{n}_1 \\ \vdots \\ \bar{n}_{m-1} \end{pmatrix}$$

$$+ [W^{2m}_{2m-1}, \dots, W^{2m}_{m+1}](\bar{z}_{m-1}, \dots, \bar{z}_0) \cdot \begin{pmatrix} \bar{z}_{2m-1} \\ \vdots \\ \bar{z}_{m+1} \end{pmatrix}$$

$$+ W^{2m}_m(\bar{z}_{m-1}, \dots, \bar{z}_0) \cdot \bar{z}_m + R^{2m}(\bar{z}_{m-1}, \dots, \bar{z}_0) + \frac{(2m)!}{2(m!)^2} G_0^2 \cdot \bar{z}_m^2$$

$$\overset{(3.4)}{=} \bar{S}_1(\bar{z}_0) \cdot n_0 + [W^{2m}_{2m-1}, \dots, W^{2m}_{m+1}](\bar{z}_{m-1}, \dots, \bar{z}_0) \cdot \begin{pmatrix} & a^{2m-1} \\ M^{2m-1} \cdot A^{2m-1} & \end{pmatrix}_{\overline{m}} \cdot \begin{pmatrix} \bar{n}_1 \\ \vdots \\ \bar{n}_{m-1} \end{pmatrix}$$

$$+ T^{2m}[0, \bar{z}_{2m-1}, \dots, \bar{z}_0]$$

$$\overset{(4.8)}{=} \bar{S}_1(\bar{z}_0) \cdot n_0 - \bar{S}_1(\bar{z}_0) \cdot \bar{z}_{2m}$$

$$+ [W^{2m-1}_{2m-2}, \dots, W^{2m-1}_m](\bar{z}_{m-1}, \dots, \bar{z}_0) \cdot \begin{pmatrix} d_{2m-1,2} & & \\ & \ddots & \\ & & d_{2m-1,m} \end{pmatrix} \cdot \begin{pmatrix} & a^{2m-1} \\ M^{2m-1} \cdot A^{2m-1} & \end{pmatrix}_{\overline{m}} \cdot \begin{pmatrix} \bar{n}_1 \\ \vdots \\ \bar{n}_{m-1} \end{pmatrix}$$

under consideration of

$$0 = T^{2m}[\bar{z}_{2m}, \dots, \bar{z}_0] \overset{(3.4)}{=} \bar{S}_1(\bar{z}_0) \cdot \bar{z}_{2m} + T^{2m}[0, \bar{z}_{2m-1}, \dots, \bar{z}_0].$$

Now by (ii) with $k$ replaced by $m-1$ and (3.20), the operator $[W^{2m-1}_{2m-2}, \dots, W^{2m-1}_m](\bar{z}_{m-1}, \dots, \bar{z}_0)$ can be substituted according to

$$0 = \bar{S}_1(\bar{z}_0) \cdot n_0 - \bar{S}_1(\bar{z}_0) \cdot \bar{z}_{2m}$$

$$+ [\bar{S}_2(\bar{z}_1, \bar{z}_0), \cdots, \bar{S}_m(\bar{z}_{m-1}, \dots, \bar{z}_0)] \cdot \begin{pmatrix} c_{2m-1,2} & & \\ & \ddots & \\ & & c_{2m-1,m} \end{pmatrix} \cdot (M^{2m-1})^{-1}$$

$$\cdot \underbrace{\begin{pmatrix} d_{2m-1,2} & & \\ & \ddots & \\ & & d_{2m-1,m} \end{pmatrix} \cdot \begin{pmatrix} & a^{2m-1} \\ M^{2m-1} \cdot A^{2m-1} & \end{pmatrix}_{\overline{m}}}_{(iv)} \cdot \begin{pmatrix} \bar{n}_1 \\ \vdots \\ \bar{n}_{m-1} \end{pmatrix}$$

and using again (iv) with $k$ replaced by $m-1$, the operator $(M^{2m-1})^{-1}$ is skipped, implying

$$0 = \bar{S}_1(\bar{z}_0)\cdot n_0 - \bar{S}_1(\bar{z}_0)\cdot \bar{z}_{2m}$$

$$+ [\,\bar{S}_2(\bar{z}_1,\bar{z}_0)\,,\cdots,\bar{S}_m(\bar{z}_{m-1},\dots,\bar{z}_0)\,]\cdot \underbrace{\begin{pmatrix} c_{2m-1,2} & & \\ & \ddots & \\ & & c_{2m-1,m} \end{pmatrix}\cdot\begin{pmatrix} d_{2m-1,2} & & \\ & \ddots & \\ & & d_{2m-1,m} \end{pmatrix}}_{(3.13)}\cdot\begin{pmatrix} \bar{n}_1 \\ \vdots \\ \bar{n}_{m-1} \end{pmatrix}$$

$$= \bar{S}_1(\bar{z}_0)\cdot n_0 - \bar{S}_1(\bar{z}_0)\cdot \bar{z}_{2m} + [\,\bar{S}_2(\bar{z}_1,\bar{z}_0)\,,\cdots,\bar{S}_m(\bar{z}_{m-1},\dots,\bar{z}_0)\,]\cdot\begin{pmatrix} c_{2m,2} & & \\ & \ddots & \\ & & c_{2m,m} \end{pmatrix}\cdot\begin{pmatrix} \bar{n}_1 \\ \vdots \\ \bar{n}_{m-1} \end{pmatrix}$$

$$\overset{(4.3)}{\underset{(3.14)}{=}} [\,\bar{S}_1(\bar{z}_0)\,,\ \bar{S}_2(\bar{z}_1,\bar{z}_0)\,,\cdots,\bar{S}_m(\bar{z}_{m-1},\dots,\bar{z}_0)\,]\cdot C^{2m}\cdot\begin{pmatrix} n_0 \\ \bar{n}_1 \\ \vdots \\ \bar{n}_{m-1} \end{pmatrix} - \bar{S}_1(\bar{z}_0)\cdot \bar{z}_{2m}\,. \qquad (4.16)$$

Hence, considering

$$N\Big[\,[\,\bar{S}_1\,,\cdots,\,\bar{S}_m\,]\cdot C^{2m}{}_{|N_0\times\cdots\times N_{m-1}}\Big] = N\big[\,[\,\bar{S}_1\,,\cdots,\,\bar{S}_m\,]_{|N_0\times\cdots\times N_{m-1}}\big]\,,$$

the affine subspace of solutions of (4.15) is given by (iii)

$$\begin{pmatrix} n_0 \\ \bar{n}_1 \\ \vdots \\ \bar{n}_{m-1} \end{pmatrix} = \begin{pmatrix} \bar{z}_{2m} \\ 0 \\ \vdots \\ 0 \end{pmatrix} + \begin{pmatrix} a^{2m} & \bar{a}^{2m} \\ A^{2m} & \bar{A}^{2m} \end{pmatrix}\cdot\begin{pmatrix} \tilde{n}_1 \\ \tilde{n}_2 \\ \vdots \\ \tilde{n}_m \end{pmatrix} \qquad (4.17)$$

with $(\tilde{n}_1,\dots,\tilde{n}_m) \in N_1\times\cdots\times N_m$ and $(\bar{z}_{2m},0,\dots,0)^T$ obviously representing a basic solution of (4.15) according to (4.16).

Summing up (4.13), (4.14) and (4.17), the solutions of the extended system

$$[\,T^{2m},T^{2m-1},T^{2m-2},\dots,T^{m+1}\,](z_{2m},z_{2m-1},\dots,z_{m+1},\ \bar{z}_m,\dots,\bar{z}_0) = 0$$

are uniquely determined by

$$\begin{pmatrix} z_{2m} \\ z_{2m-1} \\ \vdots \\ z_{m+1} \end{pmatrix} = \begin{pmatrix} \bar{z}_{2m} + (a^{2m}\ \ \bar{a}^{2m})\cdot\begin{pmatrix} \tilde{n}_1 \\ \vdots \\ \tilde{n}_m \end{pmatrix} \\ \begin{pmatrix} \bar{z}_{2m-1} \\ \vdots \\ \bar{z}_{m+1} \end{pmatrix} + \begin{pmatrix} a^{2m-1} \\ M^{2m-1}\cdot A^{2m-1} \end{pmatrix}_{\overline{m}}\cdot(A^{2m}\ \ \bar{A}^{2m})\cdot\begin{pmatrix} \tilde{n}_1 \\ \vdots \\ \tilde{n}_m \end{pmatrix} \end{pmatrix}$$

$$= \begin{pmatrix} \bar{z}_{2m} \\ \bar{z}_{2m-1} \\ \vdots \\ \bar{z}_{m+1} \end{pmatrix} + \underbrace{\begin{pmatrix} (a^{2m}\ \ \bar{a}^{2m}) \\ a^{2m-1}\cdot(A^{2m}\ \ \bar{A}^{2m}) \\ M^{2m-1}\cdot A^{2m-1}\cdot(A^{2m}\ \ \bar{A}^{2m}) \end{pmatrix}_{\overline{m+1}}}_{=M^{2m+1}\ by\ (3.18)}\cdot\begin{pmatrix} \tilde{n}_1 \\ \tilde{n}_2 \\ \vdots \\ \tilde{n}_m \end{pmatrix}$$

$$\Leftrightarrow \begin{pmatrix} z_{2m} \\ \vdots \\ z_{m+1} \end{pmatrix} = \begin{pmatrix} \bar{z}_{2m} \\ \vdots \\ \bar{z}_{m+1} \end{pmatrix} + M^{2m+1} \cdot \begin{pmatrix} \tilde{n}_1 \\ \vdots \\ \tilde{n}_m \end{pmatrix}$$

with $(\tilde{n}_1, \dots, \tilde{n}_m) \in N_1 \times \cdots \times N_m$, thus proving (i) with $k = m$.

**Proof of Lemma 4.1 (ii) :** In the next step, (ii) is accomplished by induction from $m = 1$ to $m = k$ by use of (iv) according to figure 1. In case of $m = 1$, the identity (ii) reads

$$W^3(\bar{z}_1, \bar{z}_0) \cdot \begin{pmatrix} I_B & \\ & M^3 \end{pmatrix} = [\, \bar{S}_1(\bar{z}_0)\,,\ \bar{S}_2(\bar{z}_1, \bar{z}_0)\,] \cdot C^3$$

yielding by (3.3), (3.6), (3.7), (3.11), (3.13) and (3.15)

$$[\, G_0^1\,,\ 3G_0^2\bar{z}_1\,] \cdot \begin{pmatrix} I_B & \\ & I_B \end{pmatrix} = [\, G_0^1\,,\ 2G_0^2\bar{z}_1\,] \cdot \begin{pmatrix} 1 & \\ & \frac{3}{2} \end{pmatrix}$$

which is true by inspection. Now, suppose identity (ii) with $k$ replaced by $m - 1 \in \{1, \dots, k-2\}$

$$W^{2m-1}(\bar{z}_{m-1}, \dots, \bar{z}_0) \cdot \begin{pmatrix} I_B & \\ & M^{2m-1} \end{pmatrix} = [\, \bar{S}_1(\bar{z}_0)\,, \cdots, \bar{S}_m(\bar{z}_{m-1}, \dots, \bar{z}_0)\,] \cdot C^{2m-1}\,. \tag{4.18}$$

Our aim is to show (ii) with $k = m$ or equivalently under consideration of (3.6) and (4.3)

$$[\, W_{2m}^{2m+1}, \dots, W_{m+1}^{2m+1}\,](\bar{z}_m, \dots, \bar{z}_0) \cdot M^{2m+1} \tag{4.19}$$

$$= [\, \bar{S}_2\,, \cdots,\ \bar{S}_{m+1}\,](\bar{z}_m, \dots, \bar{z}_0) \cdot \begin{pmatrix} c_{2m+1,2} & & \\ & \ddots & \\ & & c_{2m+1,m+1} \end{pmatrix}.$$

Further, with index $m_s$ and $m_z$ denoting column $m$ and row $m$ of a matrix respectively, as well as index $m|$ denoting a matrix without column $m$, the left hand side of (4.19) reads

$$\{\, [\, W_{2m}^{2m+1}, \dots, W_{m+1}^{2m+1}\,](\bar{z}_m, \dots, \bar{z}_0) \cdot M_{m|}^{2m+1}\,,\ \ [\, W_{2m}^{2m+1}, \dots, W_{m+1}^{2m+1}\,](\bar{z}_m, \dots, \bar{z}_0) \cdot M_{m_s}^{2m+1}\,\}$$

$$\overset{(3.19)}{=} \{\, [\, W_{2m}^{2m+1}, \dots, W_{m+2}^{2m+1}\,](\bar{z}_m, \dots, \bar{z}_0) \cdot M_{m|,\bar{m}}^{2m+1}\,, \tag{4.20}$$

$$[\, W_{2m}^{2m+1}, \dots, W_{m+2}^{2m+1}\,](\bar{z}_m, \dots, \bar{z}_0) \cdot M_{m_s,\bar{m}}^{2m+1} + W_{m+1}^{2m+1}(\bar{z}_m, \dots, \bar{z}_0) \cdot I_B\}$$

$$\overset{(4.7)}{\underset{(4.6)}{=}} \{\, [\, W_{2m-1}^{2m}, \dots, W_{m+1}^{2m}\,](\bar{z}_{m-1}, \dots, \bar{z}_0) \cdot \begin{pmatrix} d_{2m,2} & & \\ & \ddots & \\ & & d_{2m,m} \end{pmatrix} \cdot M_{m|,\bar{m}}^{2m+1}\,,$$

$$[\, W_{2m-1}^{2m}, \dots, W_{m+1}^{2m}\,](\bar{z}_{m-1}, \dots, \bar{z}_0) \cdot \begin{pmatrix} d_{2m,2} & & \\ & \ddots & \\ & & d_{2m,m} \end{pmatrix} \cdot M_{m_s,\bar{m}}^{2m+1}$$

$$+ \ [\, W_m^{2m}(\bar{z}_{m-1}, \dots, \bar{z}_0) + \frac{(2m)!}{(m!)^2}\, G_0^2 \cdot \bar{z}_m\,] \cdot d_{2m,m+1}\,\}\,.$$

The right hand side of (4.19) can be represented by induction hypothesis (4.18) according to

$$\{\, [\, \bar{S}_2, \cdots, \bar{S}_m\, ](\bar{z}_{m-1}, \dots, \bar{z}_0) \cdot \begin{pmatrix} c_{2m+1,2} & & \\ & \ddots & \\ & & c_{2m+1,m} \end{pmatrix}, \quad \bar{S}_{m+1}(\bar{z}_m, \dots, \bar{z}_0) \cdot c_{2m+1,m+1} \,\} \tag{4.21}$$

$$\overset{(4.18)\,(3.16)\,(3.2)}{=} \{\, [\, W^{2m-1}_{2m-2}, \dots, W^{2m-1}_m\, ](\bar{z}_{m-1}, \dots, \bar{z}_0) \cdot M^{2m-1} \cdot \begin{pmatrix} c_{2m-1,2} & & \\ & \ddots & \\ & & c_{2m-1,m} \end{pmatrix}^{-1} \begin{pmatrix} c_{2m+1,2} & & \\ & \ddots & \\ & & c_{2m+1,m} \end{pmatrix},$$

$$[\, W^{2m}_{2m-1}, \dots, W^{2m}_m\, ](\bar{z}_{m-1}, \dots, \bar{z}_0) \cdot \begin{pmatrix} & \bar{a}^{2m-1} \\ M^{2m-1} \cdot & \bar{A}^{2m-1} \end{pmatrix} \cdot c_{2m+1,m+1} + \frac{(2m)!}{(m!)^2}\, G_0^2 \bar{z}_m \cdot c_{2m+1,m+1} \,\}$$

$$\overset{(4.8)\,(3.13)}{=} \{\, [\, W^{2m}_{2m-1}, \dots, W^{2m}_{m+1}\, ](\bar{z}_{m-1}, \dots, \bar{z}_0) \cdot \begin{pmatrix} d_{2m-1,2} & & \\ & \ddots & \\ & & d_{2m-1,m} \end{pmatrix}^{-1} M^{2m-1} \cdot \begin{pmatrix} c_{2m-1,2} & & \\ & \ddots & \\ & & c_{2m-1,m} \end{pmatrix}^{-1}$$

$$\cdot \begin{pmatrix} c_{2m-1,2} & & \\ & \ddots & \\ & & c_{2m-1,m} \end{pmatrix} \cdot \begin{pmatrix} d_{2m-1,2} & & \\ & \ddots & \\ & & d_{2m-1,m} \end{pmatrix} \cdot \begin{pmatrix} d_{2m,2} & & \\ & \ddots & \\ & & d_{2m,m} \end{pmatrix},$$

$$[\, W^{2m}_{2m-1}, \dots, W^{2m}_{m+1}\, ](\bar{z}_{m-1}, \dots, \bar{z}_0) \cdot \begin{pmatrix} & \bar{a}^{2m-1} \\ M^{2m-1} \cdot & \bar{A}^{2m-1} \end{pmatrix}_{\bar{m}} \cdot \overbrace{c_{2m,m+1}}^{=1} \cdot d_{2m,m+1}$$

$$+ \,[\, W^{2m}_m(\bar{z}_{m-1}, \dots, \bar{z}_0) + \frac{(2m)!}{(m!)^2}\, G_0^2 \bar{z}_m\, ] \cdot 1 \cdot d_{2m,m+1} \,\}.$$

Comparing (4.20) and (4.21), sufficient conditions for equality are given by

$$\begin{pmatrix} d_{2m,2} & & \\ & \ddots & \\ & & d_{2m,m} \end{pmatrix} \cdot M^{2m+1}_{m|,\bar{m}}$$

$$= \begin{pmatrix} d_{2m-1,2} & & \\ & \ddots & \\ & & d_{2m-1,m} \end{pmatrix}^{-1} \cdot M^{2m-1} \cdot \begin{pmatrix} d_{2m-1,2} & & \\ & \ddots & \\ & & d_{2m-1,m} \end{pmatrix} \cdot \begin{pmatrix} d_{2m,2} & & \\ & \ddots & \\ & & d_{2m,m} \end{pmatrix}$$

as well as

$$\begin{pmatrix} d_{2m,2} & & \\ & \ddots & \\ & & d_{2m,m} \end{pmatrix} \cdot M^{2m+1}_{m_s,\bar{m}} = \begin{pmatrix} & \bar{a}^{2m-1} \\ M^{2m-1} \cdot & \bar{A}^{2m-1} \end{pmatrix}_{\bar{m}} \cdot d_{2m,m+1}$$

or combined under consideration of (iv) with $k = m - 1$ and (3.17), (3.18)

$$\begin{pmatrix} d_{2m,2} & & \\ & \ddots & \\ & & d_{2m,m} \end{pmatrix} \cdot M^{2m+1}_{\bar{m}}$$

$$= \Big[ \overbrace{\begin{pmatrix} & a^{2m-1} \\ M^{2m-1} \cdot & A^{2m-1} \end{pmatrix}_{\bar{m}} \cdot \begin{pmatrix} d_{2m,2} & & \\ & \ddots & \\ & & d_{2m,m} \end{pmatrix}}^{(iv)} , \begin{pmatrix} & \bar{a}^{2m-1} \\ M^{2m-1} \cdot & \bar{A}^{2m-1} \end{pmatrix}_{\bar{m}} \cdot d_{2m,m+1} \Big]$$

$$\overset{(3.18)}{\Leftrightarrow} \begin{pmatrix} d_{2m,2} & & \\ & \ddots & \\ & & d_{2m,m} \end{pmatrix} \cdot \left( \overbrace{\begin{pmatrix} & & (a^{2m} & \bar{a}^{2m}) \\ & a^{2m-1} \cdot & (A^{2m} & \bar{A}^{2m}) \\ M^{2m-1} \cdot & A^{2m-1} \cdot & (A^{2m} & \bar{A}^{2m}) \end{pmatrix}_{\overline{m+1}}}^{=M^{2m+1}} \right)_{\bar{m}}$$

$$= \begin{pmatrix} & (a^{2m-1} & \bar{a}^{2m-1}) \\ M^{2m-1} \cdot & (A^{2m-1} & \bar{A}^{2m-1}) \end{pmatrix}_{\bar{m}} \cdot \begin{pmatrix} d_{2m,2} & & \\ & \ddots & \\ & & d_{2m,m+1} \end{pmatrix}$$

$$\Leftrightarrow \begin{pmatrix} d_{2m,2} & & \\ & \ddots & \\ & & d_{2m,m} \end{pmatrix} \cdot \begin{pmatrix} & (a^{2m} \quad \bar{a}^{2m}) \\ \begin{pmatrix} & a^{2m-1} \\ M^{2m-1} \cdot & A^{2m-1} \end{pmatrix}_{\bar{m}} & \cdot (A^{2m} \quad \bar{A}^{2m}) \end{pmatrix}_{\bar{m}}$$

$$= \begin{pmatrix} & (a^{2m-1} & \bar{a}^{2m-1}) \\ M^{2m-1} \cdot & (A^{2m-1} & \bar{A}^{2m-1}) \end{pmatrix}_{\bar{m}} \cdot \begin{pmatrix} d_{2m,2} & & \\ & \ddots & \\ & & d_{2m,m+1} \end{pmatrix}.$$

Now, using again (iv), we end up with

$$\begin{pmatrix} d_{2m,2} & & \\ & \ddots & \\ & & d_{2m,m} \end{pmatrix} \left( \begin{matrix} & (a^{2m} \quad \bar{a}^{2m}) \\ \overbrace{\begin{pmatrix} d_{2m-1,2} & & \\ & \ddots & \\ & & d_{2m-1,m} \end{pmatrix}^{-1} \cdot M^{2m-1} \cdot \begin{pmatrix} d_{2m-1,2} & & \\ & \ddots & \\ & & d_{2m-1,m} \end{pmatrix}}^{(iv)} & \cdot (A^{2m} \quad \bar{A}^{2m}) \end{matrix} \right)_{\bar{m}}$$

$$= \begin{pmatrix} & (a^{2m-1} & \bar{a}^{2m-1}) \\ M^{2m-1} \cdot & (A^{2m-1} & \bar{A}^{2m-1}) \end{pmatrix}_{\bar{m}} \cdot \begin{pmatrix} d_{2m,2} & & \\ & \ddots & \\ & & d_{2m,m+1} \end{pmatrix}. \tag{4.22}$$

Further, (3.13) and (3.17) yield

$$\begin{pmatrix} a^{2m} & \bar{a}^{2m} \\ A^{2m} & \bar{A}^{2m} \end{pmatrix} = (C^{2m})^{-1} \cdot E^m \cdot C^{2m} \tag{4.23}$$

$$= (D^{2m-1})^{-1} \cdot (C^{2m-1})^{-1} \cdot E^m \cdot C^{2m-1} \cdot D^{2m-1} = (D^{2m-1})^{-1} \cdot \begin{pmatrix} a^{2m-1} & \bar{a}^{2m-1} \\ A^{2m-1} & \bar{A}^{2m-1} \end{pmatrix} \cdot D^{2m-1}$$

and for later use, we obtain by a similar calculation

$$\begin{pmatrix} a^{2m} & \bar{a}^{2m} \\ A^{2m} & \bar{A}^{2m} \end{pmatrix} = D^{2m} \cdot \begin{pmatrix} a^{2m+1} \\ A^{2m+1} \end{pmatrix}_{\overline{m+1}} \cdot (D^{2m})^{-1} . \tag{4.24}$$

Then by (4.23), the first row in (4.22) is equivalent to

$$d_{2m,2} \cdot \overbrace{\left(d_{2m-1,1}\right)^{-1}}^{=1} \cdot (a^{2m-1} \quad \bar{a}^{2m-1}) \cdot \overbrace{\begin{pmatrix} d_{2m-1,1} & & \\ & \ddots & \\ & & d_{2m-1,m} \end{pmatrix}}^{=D^{2m-1}}$$

$$= (a^{2m-1} \quad \bar{a}^{2m-1}) \cdot \begin{pmatrix} d_{2m,2} & & \\ & \ddots & \\ & & d_{2m,m+1} \end{pmatrix}$$

$$\Leftrightarrow \quad 0 \;=\; (a^{2m-1} \quad \bar{a}^{2m-1}) \cdot \underbrace{[\; d_{2m,2} \cdot \begin{pmatrix} d_{2m-1,1} & & \\ & \ddots & \\ & & d_{2m-1,m} \end{pmatrix} - \begin{pmatrix} d_{2m,2} & & \\ & \ddots & \\ & & d_{2m,m+1} \end{pmatrix} ]}_{=0 \; by \; (4.4)}$$

which is obviously true using (4.4) with $j = 0, \dots, m-1$. Along the same lines of reasoning, the remaining rows $2, \dots m-1$ in (4.22) can be treated with $j = 1, \dots, m-2$ according to

$$d_{2m,2+j} \cdot \left(d_{2m-1,1+j}\right)^{-1} \cdot M_{j_z}^{2m-1} \cdot \overbrace{(A^{2m-1} \quad \bar{A}^{2m-1}) \cdot \begin{pmatrix} d_{2m-1,1} & & \\ & \ddots & \\ & & d_{2m-1,m} \end{pmatrix}}^{(4.23)}$$

$$= M_{j_z}^{2m-1} \cdot (A^{2m-1} \quad \bar{A}^{2m-1}) \cdot \begin{pmatrix} d_{2m,2} & & \\ & \ddots & \\ & & d_{2m,m+1} \end{pmatrix}$$

$$\Leftrightarrow$$

$$0 = M_{j_z}^{2m-1} \cdot (A^{2m-1} \quad \bar{A}^{2m-1}) \cdot \underbrace{[\; \overbrace{\left(\frac{d_{2m,2+j}}{d_{2m-1,1+j}}\right)}^{=d_{2m,2} \; by \; (4.4)} \cdot \begin{pmatrix} d_{2m-1,1} & & \\ & \ddots & \\ & & d_{2m-1,m} \end{pmatrix} - \begin{pmatrix} d_{2m,2} & & \\ & \ddots & \\ & & d_{2m,m+1} \end{pmatrix} ]}_{=0 \; by \; (4.4)} ,$$

thus finishing the proof of Lemma 4.1 (ii).

**Proof of Lemma 4.1 (iv) :** The proof of (iv) is accomplished independently of (i)-(iii) using an inductive argument, again from $m = 1$ to $m = k$. For $m = 1$, identity (iv) is obviously true by

$$\underbrace{M^3}_{=I_B \; by \; (3.15)} \cdot d_{3,2} = d_{3,2} \cdot \underbrace{a^3}_{=I_B \; by \; (3.15)}$$

Then, suppose (iv) with $k$ replaced by $m - 1 \in \{1, \dots, k-2\}$ yielding

$$M^{2m-1}\cdot\begin{pmatrix} d_{2m-1,2} & & \\ & \ddots & \\ & & d_{2m-1,m}\end{pmatrix} = \begin{pmatrix} d_{2m-1,2} & & \\ & \ddots & \\ & & d_{2m-1,m}\end{pmatrix}\cdot\begin{pmatrix} & a^{2m-1} \\ M^{2m-1}\cdot & A^{2m-1}\end{pmatrix}_{\overline{m}}, \tag{4.25}$$

whereas the identity to prove is given by

$$\overbrace{\begin{pmatrix} & (a^{2m} \quad \bar{a}^{2m}) \\ \begin{pmatrix} & a^{2m-1} \\ M^{2m-1}\cdot & A^{2m-1}\end{pmatrix} & \cdot\,(A^{2m} \quad \bar{A}^{2m})\end{pmatrix}_{\overline{m+1}}}^{=M^{2m+1}} \cdot\begin{pmatrix} d_{2m+1,2} & & \\ & \ddots & \\ & & d_{2m+1,m+1}\end{pmatrix} \tag{4.26}$$

$$= \begin{pmatrix} d_{2m+1,2} & & \\ & \ddots & \\ & & d_{2m+1,m+1}\end{pmatrix}\cdot\begin{pmatrix} & a^{2m+1} \\ \underbrace{\begin{pmatrix} & (a^{2m} \quad \bar{a}^{2m}) \\ \begin{pmatrix} & a^{2m-1} \\ M^{2m-1}\cdot & A^{2m-1}\end{pmatrix} & \cdot\,(A^{2m} \quad \bar{A}^{2m})\end{pmatrix}_{\overline{m+1}}}_{=M^{2m+1}} & \cdot\,A^{2m+1}\end{pmatrix}_{\overline{m+1}}$$

under consideration of (3.18). Further, using the induction hypothesis (4.25), the right hand side of (4.26) transforms according to

$$\begin{pmatrix} d_{2m+1,2} & & \\ & \ddots & \\ & & d_{2m+1,m+1}\end{pmatrix} \tag{4.27}$$

$$\cdot\begin{pmatrix} & a^{2m+1} \\ & (a^{2m} \quad \bar{a}^{2m})\cdot A^{2m+1} \\ \begin{pmatrix} d_{2m-1,2} & & \\ & \ddots & \\ & & d_{2m-1,m}\end{pmatrix}^{-1}\cdot\, M^{2m-1}\cdot\begin{pmatrix} d_{2m-1,2} & & \\ & \ddots & \\ & & d_{2m-1,m}\end{pmatrix} & \cdot\,(A^{2m} \quad \bar{A}^{2m})\cdot A^{2m+1}\end{pmatrix}_{\overline{m+1}}.$$

Next, the left hand side in (4.26) is transformed into (4.27), row by row. Concerning the first row, we obtain from the left hand side in (4.26) and the second relation in (4.4)

$$(a^{2m} \quad \bar{a}^{2m})\cdot\begin{pmatrix} d_{2m+1,2} & & \\ & \ddots & \\ & & d_{2m+1,m+1}\end{pmatrix}$$

$$\overset{(4.24)}{=} d_{2m,1}\cdot a^{2m+1}\cdot\begin{pmatrix} d_{2m,1} & & \\ & \ddots & \\ & & d_{2m,m}\end{pmatrix}^{-1}\cdot\begin{pmatrix} d_{2m+1,2} & & \\ & \ddots & \\ & & d_{2m+1,m+1}\end{pmatrix}$$

$$\overset{(4.4)}{=} \underbrace{d_{2m,1}}_{=1} \cdot a^{2m+1} \cdot \begin{pmatrix} d_{2m,1} & & \\ & \ddots & \\ & & d_{2m,m} \end{pmatrix}^{-1} \cdot \begin{pmatrix} d_{2m,1} & & \\ & \ddots & \\ & & d_{2m,m} \end{pmatrix} \cdot d_{2m+1,2}$$

$$= d_{2m+1,2} \cdot a^{2m+1},$$

obviously agreeing with the first row in (4.27). The second row from the left hand side in (4.26) implies

$$a^{2m-1} \cdot (A^{2m} \quad \bar{A}^{2m}) \cdot \begin{pmatrix} d_{2m+1,2} & & \\ & \ddots & \\ & & d_{2m+1,m+1} \end{pmatrix}$$

$$\overset{(4.24)}{\underset{(4.23)}{=}} \underbrace{d_{2m-1,1}}_{=1} \cdot a^{2m} \cdot \begin{pmatrix} d_{2m-1,1} & & \\ & \ddots & \\ & & d_{2m-1,m-1} \end{pmatrix}^{-1} \cdot \begin{pmatrix} d_{2m,2} & & \\ & \ddots & \\ & & d_{2m,m} \end{pmatrix} \cdot A_{\overline{m}}^{2m+1} \cdot \begin{pmatrix} d_{2m,1} & & \\ & \ddots & \\ & & d_{2m,m} \end{pmatrix}^{-1}$$

$$\cdot \begin{pmatrix} d_{2m+1,2} & & \\ & \ddots & \\ & & d_{2m+1,m+1} \end{pmatrix}$$

$$\overset{(4.4)}{=} a^{2m} \cdot d_{2m,2} \cdot A_{\overline{m}}^{2m+1} \cdot d_{2m+1,2}$$

$$\overset{(4.5)}{=} d_{2m+1,3} \cdot a^{2m} \cdot A_{\overline{m}}^{2m+1}$$

$$= d_{2m+1,3} \cdot (a^{2m} \quad \bar{a}^{2m}) \cdot A^{2m+1}$$

yielding the second row in (4.27). Note that the last identity follows from the fact that the $m$-th row of $A^{2m+1}$ equals zero according to the definition in (3.17).

Finally, in case of $m \geq 3$, the remaining rows $3, \dots m$ in (4.26), (4.27) can be treated along the same lines of reasoning with $j = 1, \dots, m-2$ according to

$$M_{j_z}^{2m-1} \cdot A^{2m-1} \cdot (A^{2m} \quad \bar{A}^{2m}) \cdot \begin{pmatrix} d_{2m+1,2} & & \\ & \ddots & \\ & & d_{2m+1,m+1} \end{pmatrix}$$

$$\overset{(4.24)}{\underset{(4.23)}{=}} M_{j_z}^{2m-1} \cdot \begin{pmatrix} d_{2m-1,2} & & \\ & \ddots & \\ & & d_{2m-1,m} \end{pmatrix} \cdot A^{2m} \cdot \begin{pmatrix} d_{2m-1,1} & & \\ & \ddots & \\ & & d_{2m-1,m-1} \end{pmatrix}^{-1}$$

$$\cdot \begin{pmatrix} d_{2m,2} & & \\ & \ddots & \\ & & d_{2m,m} \end{pmatrix} \cdot A_{\overline{m}}^{2m+1} \cdot \begin{pmatrix} d_{2m,1} & & \\ & \ddots & \\ & & d_{2m,m} \end{pmatrix}^{-1} \cdot \begin{pmatrix} d_{2m+1,2} & & \\ & \ddots & \\ & & d_{2m+1,m+1} \end{pmatrix}$$

$$\overset{(4.4)}{=} M_{j_z}^{2m-1} \cdot \begin{pmatrix} d_{2m-1,2} & & \\ & \ddots & \\ & & d_{2m-1,m} \end{pmatrix} \cdot A^{2m} \cdot d_{2m,2} \cdot A_{\overline{m}}^{2m+1} \cdot d_{2m+1,2}$$

$$\overset{(4.5)}{=} \frac{d_{2m+1,3+j}}{d_{2m-1,1+j}} \cdot M_{j_z}^{2m-1} \cdot \begin{pmatrix} d_{2m-1,2} & & \\ & \ddots & \\ & & d_{2m-1,m} \end{pmatrix} \cdot A^{2m} \cdot A_{\overline{m}}^{2m+1}$$

$$= \frac{d_{2m+1,3+j}}{d_{2m-1,1+j}} \cdot M_{j_z}^{2m-1} \cdot \begin{pmatrix} d_{2m-1,2} & & \\ & \ddots & \\ & & d_{2m-1,m} \end{pmatrix} \cdot \left(A^{2m} \quad \bar{A}^{2m}\right) \cdot A^{2m+1}$$

finishing the proof of Lemma 4.1.

*Acknowledgement*

The paper was presented on occasion of the seminar Bifurcation Theory of the ETH Zürich in Valbella, 1996. Special thanks to U. Kirchgraber, W.-J. Beyn and in particular E. Bohl.

*Matthias Stiefenhofer*

*University of Applied Sciences*

*87435 Kempten (Germany)*

*matthias.stiefenhofer@hs-kempten.de*